\documentclass{article}

\usepackage{amsfonts}
\usepackage{amssymb}
\usepackage{enumerate}
\usepackage[T1]{fontenc}

\newtheorem{ttt}{Theorem}[section]
\newtheorem{llll}[ttt]{Lemma}
\newtheorem{ccc}[ttt]{Claim}
\newtheorem{eee}[ttt]{Example}
\newtheorem{ppp}[ttt]{Proposition}
\newtheorem{rrr}[ttt]{Remark}
\newtheorem{sss}[ttt]{Statement}
\newtheorem{ddd}[ttt]{Definition}
\newtheorem{qqq}[ttt]{Question}
\newtheorem{cccc}[ttt]{Corollary}
\newtheorem{fff}[ttt]{Fact}

\newcommand{\bt}{\begin{ttt}}
\newcommand{\bl}{\begin{llll}}
\newcommand{\bcl}{\begin{ccc}}
\newcommand{\bex}{\begin{eee}}
\newcommand{\bprop}{\begin{ppp}}
\newcommand{\br}{\begin{rrr}\upshape}
\newcommand{\bs}{\begin{sss}}
\newcommand{\bd}{\begin{ddd}\upshape}
\newcommand{\bq}{\begin{qqq}}
\newcommand{\bcor}{\begin{cccc}}
\newcommand{\bfac}{\begin{fff}}

\newcommand{\bp}{\noindent\textbf{Proof. }}

\newcommand{\et}{\end{ttt}}
\newcommand{\el}{\end{llll}}
\newcommand{\ecl}{\end{ccc}}
\newcommand{\eex}{\end{eee}}
\newcommand{\epr}{\end{ppp}}
\newcommand{\er}{\end{rrr}}
\newcommand{\es}{\end{sss}}
\newcommand{\ed}{\end{ddd}}
\newcommand{\eq}{\end{qqq}}
\newcommand{\ecor}{\end{cccc}}
\newcommand{\efac}{\end{fff}}

\newcommand{\ep}{\hspace{\stretch{1}}$\square$\medskip}
\newcommand{\beq}{\begin{equation}}
\newcommand{\eeq}{\end{equation}}

\newcommand{\lab}[1]{\label{#1}}

\newcommand{\NN}{\mathbb{N}}

\newcommand{\RR}{\mathbb{R}}

\newcommand{\al}{\alpha}

\newcommand{\om}{\omega}
\newcommand{\si}{\sigma}

\newcommand{\iA}{\mathcal{A}}
\newcommand{\iB}{\mathcal{B}}

\newcommand{\iH}{\mathcal{H}}
\newcommand{\iM}{\mathcal{M}}

\newcommand{\iI}{\mathcal{I}}

\newcommand{\iN}{\mathcal{N}}

\newcommand{\C}{C[0,1]}
\newcommand{\Rn}{$\RR^n$}

\newcommand{\Hd}{\mathcal{H}^d}

\title{Hausdorff measures of different dimensions are isomorphic under 
the Continuum Hypothesis}

\begin{document}

\author{M\'arton Elekes\thanks{Partially supported by Hungarian Scientific
Foundation grant no. 37758 and F 43620.}}

\maketitle 

\begin{abstract}
We show that the Continuum Hypothesis implies that 
for every $0<d_1\leq d_2<n$ the measure spaces
$\left(\RR^n,\iM_{\iH^{d_1}},\iH^{d_1}\right)$ and
$\left(\RR^n,\iM_{\iH^{d_2}},\iH^{d_2}\right)$ are isomorphic, where
$\iH^d$ is $d$-dimensional Hausdorff measure and $\iM_{\Hd}$ is the
$\sigma$-algebra of measurable sets with respect to $\Hd$.  This is
motivated by the well-known question (circulated by D.~Preiss) whether
such an isomorphism exists if we replace measurable sets by Borel
sets.

We also investigate the related question whether every continuous
function (or the typical continuous function) is H\"older continuous
(or is of bounded variation) on a set of positive Hausdorff dimension.
\end{abstract}

\insert\footins{\footnotesize{MSC codes: Primary 28A78; Secondary
26A16, 26A45, 28A05, 03E50}} \insert\footins{\footnotesize{Key Words:
Hausdorff measure, isomorphism, $CH$, Hausdorff dimension, 
H\"older continuous, bounded variation}}

\section*{Introduction}

The following problem is circulated by D.~Preiss \cite{Pr} 
(while it is unclear, 
who actually asked this first, see also \cite{Cs}, where the question 
is under the name of Preiss).
Let $\Hd$ denote $d$-dimensional Hausdorff measure, see
e.g.~\cite{Fa}, \cite{Fe} or \cite{Ma}, and let $\iB$ denote the
$\sigma$-algebra of Borel subsets of $\RR^n$. By \emph{isomorphism} 
of two measure spaces we mean a bijection $f$ such that
both $f$ and $f^{-1}$ are measurable set and measure preserving. 

\bq\lab{David}
Let $0<d_1 < d_2<n$. Are the measure spaces
$\left(\RR^n,\iB,\iH^{d_1}\right)$ and
$\left(\RR^n,\iB,\iH^{d_2}\right)$ isomorphic?
\eq
 
An equally natural question is whether such an isomorphism exists if
we replace Borel sets by measurable sets with respect to Hausdorff
measures. Denote by $\iM_{\Hd}$ the $\sigma$-algebra of measurable
sets with respect to $\Hd$, in the usual sense of Carath\'eodory.
(Again, and throughout the paper we follow standard terminology that 
can be found e.g.~in \cite{Fa}, \cite{Fe} or \cite{Ma}.) 

\bq\lab{modified}
Let $0<d_1 < d_2<n$. Are the measure spaces
$\left(\RR^n,\iM_{\iH^{d_1}},\iH^{d_1}\right)$ and
$\left(\RR^n,\iM_{\iH^{d_2}},\iH^{d_2}\right)$ isomorphic?
\eq

The main result of the first part of this 
paper is the following affirmative answer to this question assuming $CH$; that
is, under the Continuum Hypothesis. The original Question 
\ref{David} remains open.

\bt\lab{main}
Under the Continuum Hypothesis, for every $0<d_1\leq d_2<n$ the measure spaces
$\left(\RR^n,\iM_{\iH^{d_1}},\iH^{d_1}\right)$ and
$\left(\RR^n,\iM_{\iH^{d_2}},\iH^{d_2}\right)$ are isomorphic.
\et

We remark here that in the definition of isomorphism 
by measure preserving we could mean two different
things, as we can either require that the measure of every measurable
set is preserved, or that the outer measure of every set is
preserved. However, as Hausdorff measures are Borel regular (and Borel
sets are measurable with respect to Hausdorff measures), there is no
possibility of confusion here, since measure preserving for measurable
sets is equivalent to outer measure preserving for all sets. Note also 
that the term ``Borel isomorphism'' will refer to a completely 
different notion.

We do not know if the assumption of the Continuum Hypothesis can be
dropped in Theorem \ref{main}. However, it is rather unlikely as 
the following remarks show. 
The result resembles Erd\H os-Sierpi\'nski Duality (see
e.g.~\cite{Ox}); that is, the existence under $CH$ of a bijection
which maps Lebesgue nullsets onto sets of first category and vice
versa. Denote $\iN$ and $\iM$ the $\sigma$-ideals of Lebesgue nullsets
and sets of first category, respectively. For a $\sigma$-ideal $\iI$,
the symbol $\textrm{non}(\iI)$ stands for the smallest cardinal
$\kappa$ for which there exists a set of cardinality $\kappa$ that is
not in $\iI$. Now the fact that the role of $CH$ is essential for 
Erd\H os-Sierpi\'nski Duality
follows from the well-know result that in some models of set theory
$\textrm{non}(\iN) \neq \textrm{non}(\iM)$ holds (see e.g.~\cite{BJ}).

In our present setting, denote 
$\iN_{\iH^d}$ the $\sigma$-ideal of negligible sets with
respect to $\iH^d$. If it were known that in some model of set theory
$\textrm{non}(\iN_{\iH^{d_1}}) \neq \textrm{non}(\iN_{\iH^{d_2}})$
held, it would be proven that the Continuum Hypothesis cannot be
dropped in Theorem \ref{main}. However, this is a very recent research
area, and so far the above statement is only known to hold in some
model when $d_2=n$ (see \cite{SS}).

We also prove in this paper that under $CH$ for every $0<d_1\leq d_2<n$ the
measure spaces $(\RR^n,\iB^{\si f}_{d_1},\iH^{d_1})$ and
$(\RR^n,\iB^{\si f}_{d_2},\iH^{d_2})$ are isomorphic, where
$\iB^{\si f}_d$ stands for the $\sigma$-algebra of Borel subsets of
$\RR^n$ that are of $\sigma$-finite or co-$\sigma$-finite
$\Hd$-measure. This result shows that Question \ref{David} 
cannot be solved by the natural approach, as it is not enough to
consider only $d_1$-sets and $d_2$-sets; that is, sets of positive and finite 
measure with respect to $\iH^{d_1}$ and $\iH^{d_2}$. 

In the second part of the paper, motivated by Question \ref{David}, we
consider the following set of problems. Can we find for every $f:[0,1]\to\RR$
continuous/Borel/typical continuous (in the Baire category sense, 
see e.g.~\cite{Br})
function a set of positive Hausdorff dimension on which the function
agrees with a function of bounded variation/Lipschitz/H\"older
continuous function? For example it is clear that showing that every Borel
function is H\"older continuous of some suitable exponent on a set of
sufficiently large Hausdorff dimension would answer Question
\ref{David} in the negative. The other versions are less closely related 
to our problem, however, they are of independent interest. 
We prove the following two results.

\bt
Fix $0<\al\leq 1$. A typical continuous function is not H\"older continuous of
exponent $\al$ on any set of Hausdorff dimension larger than $1-\al$.
\et

We do not know whether $1-\al$ is sharp.

\bt
A typical continuous function does not agree with any function of 
bounded variation on any set of Hausdorff dimension larger than $\frac{1}{2}$.
\et
 
This theorem is motivated by an analogous result of Humke and
Laczkovich \cite{HL}, who proved that a typical continuous function is
not monotonic on any set of positive Hausdorff dimension. So one would
expect the same for functions of bounded variation, and it is a very
natural problem whether this dimension is indeed $0$, or
$\frac{1}{2}$ or perhaps something in between.

\bq
Is it true that a typical continuous function agrees with a function of 
bounded variation on a set of Hausdorff dimension $\frac{1}{2}$?
\eq

\section{Preliminaries}

In this section we present the lemmas we use in the sequel. 
Probably most of them are well-known, however, we could not find suitable 
references, so we could not avoid including them. Let
$\lambda$ denote one-dimensional Lebesgue measure, then $\iM_\lambda$ is 
the class of Lebesgue measurable sets. A \emph{Borel isomorphism} is a
bijection such that images and preimages of Borel sets are Borel.

\bl\lab{isom} 
Let $d>0$, $B\subset\RR^n$ be Borel such that
$0<\Hd(B)<\infty$ and let $I=[0,\Hd(B)]$. Then there exists an isomorphism
$f$ between the measure spaces $(B,\iM_{\Hd},\Hd)$ and
$(I,\iM_\lambda,\lambda)$ that is also a Borel isomorphism.
\el

\bp 
We may clearly assume $\Hd(B)=1$.
It is stated in \cite[12.B]{Ke} that every Borel set, 
hence $B$ is a standard Borel space (that is,
$B$ is Borel isomorphic to some Polish space). Hence we can apply
\cite[17.41]{Ke} which states that every continuous (that is,
singletons are of measure zero) probability measure on a standard Borel space 
is isomorphic by a Borel isomorphism to Lebesgue measure on the unit interval.
\ep

\bl\lab{sisom} 
Let $B\subset\RR^n$ be a Borel set of infinite but $\sigma$-finite
$\Hd$-measure. Then there exists an isomorphism
$f$ between the measure spaces $(B,\iM_{\Hd},\Hd)$ and
$([0,\infty),\iM_\lambda,\lambda)$ that is also a Borel isomorphism.
\el

\bp
Define a system $\iA$ of pairwise disjoint Borel subsets of $B$ by
transfinite recursion, such that $\Hd(A)=1$ for every $A\in\iA$. By
$\sigma$-finiteness this procedure stops at some countable ordinal, so
$\iA$ is countable. Put $m=\Hd(B\setminus\cup\iA)$. As we can always 
find a Borel subset of $\Hd$-measure $1$ inside a Borel set of 
$\Hd$-measure at least $1$ (see e.g.~\cite[8.20]{Ma}), we obtain 
$m<1$. In particular, $\iA$ is infinite, say $\iA=\{A_0, A_1,\dots\}$. 
By the
previous lemma $B\setminus\cup_{i=0}^\infty A_i$ is isomorphic to the
interval $[-m,0)$ (or more precisely to $[-m,0]$, but these two
intervals are isomorphic again by the previous lemma), and $A_i$ is
isomorphic to $[i,i+1)$ for every $i$. Therefore $B$ is clearly
isomorphic to $[-m,\infty)$ which is obviously isomorphic to
$[0,\infty)$.
\ep

\bl\lab{Hmeas} 
Let $n\in\NN$, $d\geq 0$ and $H\subset\RR^n$ be
arbitrary. Then the following statements are equivalent:

\begin{enumerate}[(i)]
\item\lab{egy} $H$ is $\Hd$-measurable,
\item\lab{ketto} $\Hd (B)=\Hd (B\cap H)+\Hd (B\cap H^C)$ for every
Borel set $B\subset\RR^n$ with $0<\Hd (B)<\infty$,
\item\lab{harom} $H\cap B$ is $\Hd$-measurable for every Borel
set $B\subset\RR^n$ with $0<\Hd (B)<\infty$,
\item\lab{uccso} For every Borel set $B\subset\RR^n$ with $0<\Hd (B)<\infty$, we
have $H\cap B=A\cup N$, where $A$ is Borel and $N$ is
$\Hd$-negligible.
\end{enumerate}
\el

\bp (\ref{egy}) $\iff$ (\ref{ketto}) By definition, the set $H$ is 
$\Hd$-measurable in the sense of
Carath\'eodory if and only if $\Hd(X)=\Hd(X\cap H)+\Hd(X\cap H^C)$ for every
$X\subset \RR^n$. As outer measures are subadditive, this is  equivalent to
$\Hd(X)\geq\Hd(X\cap H)+\Hd(X\cap H^C)$ for every
$X\subset \RR^n$. Once this inequality fails to hold, using the Borel
regularity of Hausdorff measures (see e.g.~\cite[4.5]{Ma}), there is a Borel set
$B\supset X$ such that $\Hd(B)<\Hd(X\cap H)+\Hd(X\cap H^C)$. Therefore $B$ is
of finite measure, moreover,
$\Hd(B)<\Hd(B\cap H)+\Hd(B\cap H^C)$, in particular, $B$ is clearly
of positive measure.

(\ref{ketto}) $\iff$ (\ref{harom}) The third condition obviously implies the
second one as $\Hd$ is additive on measurable sets. So suppose (\ref{ketto}). 
It is enough to show that $\Hd (A)\geq
\Hd(A\cap (H\cap B))+\Hd (A\cap (H\cap B)^C)$ for every
Borel set $A\subset\RR^n$ with $0<\Hd (A)<\infty$. We can assume that
$\Hd(A\cap B)>0$, otherwise the above inequality clearly holds, as 
the first term on
the right hand side is $0$, while the second term is not greater than the left
hand side. Now
\[
\Hd(A) \geq
\Hd(A\cap B) + \Hd(A\cap B^C) \geq
\]
\[
\Hd((A\cap B)\cap H) + \Hd((A\cap B)\cap H^C) + \Hd(A\cap B^C) \geq
\]
\[
\Hd(A\cap (H\cap B)) + \Hd(A\cap (H\cap B)^C),
\]
where the first inequality holds as $B$ is measurable, the second one is
(\ref{ketto}) applied to $A\cap B$, and the third one follows from the
subadditivity of $\Hd$ since $A\cap (H\cap B)^C$ is the (disjoint) union of
$(A\cap B)\cap H^C$ and $A\cap B^C$.

(\ref{harom}) $\iff$ (\ref{uccso}) As negligible sets are measurable, the
fourth condition implies the third, while we immediately obtain the other
direction if we apply Lemma \ref{isom} to $B$.
\ep

\br 
We could also assume that the above $B$ is compact, but we
will not need this fact.
\er



The following lemma is essentially \cite[2.5.10]{Fe}.

\bl\lab{loc} 
Let $0<d<n$ and suppose $CH$ holds. Then there
exists a disjoint family $\{B_\al:\al<\om_1\}$ of Borel subsets
of $\RR^n$ of finite $\Hd$-measure, such that a set $H\subset\RR^n$ is
$\Hd$-measurable iff $H\cap B_\al$ is $\Hd$-measurable for every
$\al<\om_1$.  
\el

\bp 
Let $\{A_\al:\al<\om_1\}$ be an enumeration of the Borel subsets of
\Rn\ of positive finite $\Hd$-measure, and put
$B_\al=A_\al\setminus(\cup_{\beta<\al}A_\beta)$. These are clearly pairwise
disjoint Borel sets of finite $\Hd$-measure. 
The other
direction being trivial we only have to verify that if $H\subset\RR^n$ is
such that $H\cap B_\al$ is $\Hd$-measurable for every $\al<\om_1$, then $H$
itself is $\Hd$-measurable. By Lemma \ref{Hmeas} we only have to show that
$H\cap A_\alpha$ is $\Hd$-measurable for every $\al<\om_1$. But
$A_\al=\cup_{\beta\le\al}B_\beta$, therefore $H\cap
A_\al=\cup_{\beta\le\al}(H\cap B_\al)$, which is clearly
$\Hd$-measurable, which completes the proof.
\ep

\bl\lab{sloc} 
Let $0<d<n$ and suppose $CH$ holds. Then there exists a disjoint family
$\{S_\al:\al<\om_1\}$ of Borel subsets of $\RR^n$ of infinite but
$\sigma$-finite $\Hd$-measure, such that a set $H\subset\RR^n$ is
$\Hd$-measurable iff $H\cap S_\al$ is $\Hd$-measurable for every $\al<\om_1$.
\el

\bp First we check that uncountably many
$B_\alpha$ of Lemma \ref{loc} are of positive $\Hd$-measure. 
Otherwise, as $\RR^n$ is not $\si$-finite and as by \cite[8.20]{Ma} 
every Borel set of infinite $\Hd$-measure contains a Borel set of 
$\Hd$-measure $1$, we could find a
Borel set of positive and finite measure that is disjoint from these
$B_\alpha$. But this set was enumerated as $A_\al$ for some $\alpha$,
moreover $A_\al=\cup_{\beta\le\al} B_\beta$. Since $A_\al$ is disjoint
from all $B_\beta$ of positive measure, by $CH$, countably
many zerosets cover $A_\al$, so it is a zeroset, a contradiction. 

This obviously implies that for some integer $N$ uncountably many
$B_\alpha$ are of measure at least $\frac{1}{N}$. Now we can
recursively define a partition of the set of countable ordinals into
countable intervals $\{I_\al:\al<\om_1\}$ such that for every $\al$
the $\Hd$-measure of $\cup_{\xi\in I_\al} B_\xi$ is infinite. On the
other hand, this measure is clearly $\si$-finite. Now we check that
$S_\al = \cup_{\xi\in I_\al} B_\xi$ works. Every $S_\al$ is clearly a Borel
subsets of $\RR^n$ of infinite but $\sigma$-finite $\Hd$-measure. Now
we have to show that the last statement of the lemma holds, namely
that $\Hd$-measurability ``reflects'' on the sets $S_\al$. One
direction is trivial, so in order to prove the other one let us assume
that $H\cap S_\al$ is $\Hd$-measurable for every $\al<\om_1$. We have
to show that $H$ is $\Hd$-measurable. But this is obvious by the
previous lemma, as every $B_\al$ can be covered by some $S_\beta$,
hence $H\cap B_\al = B_\al \cap (H\cap S_\beta)$ which is measurable. 
\ep

\bl\lab{sfinite}
Assume $CH$. A set $H\subset \RR^n$ is of $\sigma$-finite $\Hd$-measure iff it
can be covered by countably many of the above $S_\al$.
\el

\bp
One direction is trivial. For the other one we can assume that $H\subset
\RR^n$ is of finite $\Hd$-measure. There exists a Borel set $B$ of positive 
finite
$\Hd$-measure containing $H$. This set was enumerated in Lemma \ref{loc} as
$A_\al$ for some $\al$, so by $CH$ it can be covered by countably many $B_\al$
hence also by countably many $S_\al$.
\ep

\section{Isomorphic measure spaces}\lab{s:isom}

In this section we use the above lemmas to prove the existence of 
isomorphisms under the Continuum Hypothesis.

\bt
Under $CH$ for every $0<d_1\leq d_2<n$ the measure spaces
$\left(\RR^n,\iM_{\iH^{d_1}},\iH^{d_1}\right)$ and
$\left(\RR^n,\iM_{\iH^{d_2}},\iH^{d_2}\right)$ are isomorphic.
\et

\bp
Find two partitions $\{S_\al^{d_1}:\al<\om_1\}$ and $\{S_\al^{d_2}:\al<\om_1\}$
of $\RR^n$ as in Lemma \ref{sloc}. By Lemma \ref{sisom} find isomorphisms
$f_\al:S_\al^{d_1} \to S_\al^{d_2}$ for every $\al$. Define $f =
\cup_{\al<\om_1} f_\al$. We have to check that $f$ is an isomorphism. First we prove that $f$ preserves
measurable sets. Suppose $H\in \iM_{\iH^{d_1}}$. By
Lemma \ref{sloc} it is sufficient to show that $f(H)\cap S_\al^{d_2} \in
\iM_{\iH^{d_2}}$ for every $\al$. But $f(H)\cap S_\al^{d_2} = f_\al(H\cap
S_\al^{d_1})$ which is $\iH^{d_2}$-measurable as $f_\al$ is an
isomorphism. Similarly, $f^{-1}$ also preserves measurable sets. Now we show that
$f$ preserves measure. (Again, the same argument works
for $f^{-1}$.) As mentioned in the Introduction it is enough to show this for
measurable sets. First, 
$f$ preserves non-$\si$-finiteness since it clearly preserves
the property that characterizes $\si$-finiteness in Lemma
\ref{sfinite}, so we can restrict ourselves
to $\si$-finite sets. But such a set is partitioned by the countably many
$S_\al^{d_1}$ that cover it, and the countably many isomorphisms $f_\al$
preserve measure, so the proof is complete.
\ep

\br If $d=0$ then all subsets of $\RR^n$ are measurable, while if $d=n$ then
$\iH^d$ is $\si$-finite, therefore the theorem cannot be extended to these
cases.

As we already mentioned in the Introduction, 
it is unknown whether $CH$ can be dropped from the theorem, but 
the paper \cite{SS} is a huge step in this direction.
\er

Now we prove the existence of another type of isomorphism. Let $\iB^{\si f}_d$ 
denote the $\sigma$-algebra of Borel subsets of
$\RR^n$ that are of $\sigma$-finite or co-$\sigma$-finite $\Hd$-measure.

\bt
Under $CH$
for every $0<d_1\leq d_2<n$ the measure spaces
$\left(\RR^n,\iB^{\si f}_{d_1},\iH^{d_1}\right)$ and
$\left(\RR^n,\iB^{\si f}_{d_2},\iH^{d_2}\right)$ are isomorphic.
\et

\bp Let $f$ be as above.
It is sufficient to show that $B\in\iB^{\si f}_{d_1}$ implies $f(B)\in\iB^{\si f}_{d_2}$.
It is also sufficient to show this for $\si$-finite sets, as then for $B$
co-$\si$-finite we obtain $f(\RR^n\setminus B)\in\iB^{\si f}_{d_2}$ which
implies $f(B) = \RR^n\setminus f(\RR^n\setminus B) \in\iB^{\si f}_{d_2}$.
So let $B$ be a $\si$-finite Borel set. By Lemma
\ref{sfinite} the set $B$ is partitioned by the countably many
Borel sets $S_\al^{d_1}$ that cover it, and the countably many Borel
isomorphisms $f_\al$ produce a $\si$-finite Borel image.
\ep

This result shows that Question \ref{David} 
cannot be solved by the natural approach, as it is not enough to
consider only $d_1$-sets and $d_2$-sets; that is, sets of positive and finite 
measure with respect to $\iH^{d_1}$ and $\iH^{d_2}$.

\section{Typical continuous functions}\lab{typical}

In this section we consider a set of problems related to the previous
section. The reason is that a possible way to settle the isomorphism
problem might be answering the following question. For the sake of
simplicity, in this section we restrict ourselves to functions defined
on $[0,1]$ instead of $\RR^n$.

\bq 
For which pair of reals $0<\al,\beta\leq 1$ is it true that for
every Borel function $f:[0,1]\to \RR$ there exists a set $H\subset [0,1]$ of
Hausdorff dimension at least $\beta$ such that $f$ restricted to $H$
is H\"older continuous of exponent $\al$?
\eq

(Again, see \cite{Fa}, \cite{Fe} or \cite{Ma} for definitions). As every Borel function is continuous on a set of positive measure, we can ask
the same question for continuous functions as well. Here typical continuous
functions (in the Baire category sense, see e.g.~\cite{Br}) are natural candidates for badly
behaving examples. Indeed, the following holds.

\bt\lab{Holder}
Fix $0<\al\leq 1$. A typical continuous function is not H\"older continuous of
exponent $\al$ on any set of Hausdorff dimension larger than $1-\al$.
\et

\bp
Denote $\C$ the Banach space
 of continuous real-valued functions defined on $[0,1]$.
We say that $g\in\C$ is H\"older continuous of exponent $\al$ and constant $K$
if $|g(x)-g(y)| \leq K|x-y|^\al$ for every $x,y\in [0,1]$. 
Length of an interval $I$ is denoted by $|I|$. For $f,g\in\C$, the set $\{x\in[0,1] : f(x)=g(x)\}$ is abbreviated by $\{f=g\}$.
For $H\subset[0,1]$ denote
\[
\iH^d_\infty (H)=\inf\left\{\sum_{n=1}^\infty |I_n|^d : 
\{I_n\}_{n=1}^\infty 
\textrm{ is a sequence of intervals}, H\subset\cup_{n=1}^\infty I_n\right\}.
\]
So this is the usual Hausdorff measure with no restriction on the lengths
of the covering intervals, sometimes called Hausdorff capacity. 
It is well-known and easy to check that $\iH^d(H)=0$ if and only 
if $\iH^d_\infty(H)=0$.

First we show that it is sufficient to prove that for all integers 
$N,M>0$ the set
\[
D(N,M) = \{f\in\C : \iH^{1-\al+\frac{1}{N}}_\infty (\{f=g\})<\frac{1}{M} 
\]
\[
\textrm{ for every H\"older function } g \textrm{ of exponent } \al \textrm{ and constant } 1 \}
\]
contains a dense open subset of $\C$. Indeed, this implies that the set 
$D=\cap_{N=1}^\infty \cap_{M=1}^\infty D(N,M)$ contains a dense
$G_\delta$ set, hence is residual, which in turn implies that the typical
continuous function does not agree with a H\"older continuous function
of exponent $\al$ and constant $1$ on any set of Hausdorff dimension
larger than $1-\al$. Moreover, the map $f\mapsto Kf$ is a
homeomorphism of $\C$, hence $KD=\{Kf : f\in D\}$ is residual for every $K$. 
Therefore $\cap_{K=1}^\infty KD$ is also residual, and so the typical
continuous function does not agree with a H\"older continuous function
of exponent $\al$ on any set of Hausdorff dimension
larger than $1-\al$.

Now we show that we can assume $0<\al<1$. Indeed, if the statement of 
the theorem holds for $1-\frac{1}{L}$ for every $L$, then intersecting 
the corresponding sequence of residual subsets of $\C$ we obtain the 
case $\al=1$. 

Now what remains to be proven is 
that $D(N,M)$ contains a dense open set. Let $f_0\in\C$
and $r_0>0$ be given. The closed ball centered at $f_0$ and of radius $r_0$ is
denoted by $B(f_0,r_0)$. We have to find $f_1\in\C$ and $r_1>0$ such
that $B(f_1,r_1) \subset B(f_0,r_0) \cap D(N,M)$. By uniform
continuity of $f_0$, for a large enough integer $m$, the inequality
$|x-y|<\frac{2}{m}$ implies $|f_0(x)-f_0(y)|<\frac{r_0}{5}$. The exact
value of $m$ will be chosen later. Let $k$ be another positive integer to be
fixed later. (For those who like to see the explicit choice of the 
constants in advance, $k$ and $m$ will be chosen so that, in addition to 
the above requirement concerning uniform continuity, the inequalities
(\ref{gyok2}) and (\ref{pre}) below are satisfied. The straightforward 
calculation that this choice is indeed possible can be found in the 
last paragraph of the proof.)

Now we define the piecewise constant function $\bar{f_1}$ as
follows.
For a pair of integers $0\leq i <m$ and $0\leq j <k$ set
\[
\bar{f_1}(x) = f_0(\frac{i}{m}) + j\frac{r_0}{5k} \textrm{ \ \ \ \ \  \ 
for\ \ \ \ \ \ }
\frac{i}{m}+\frac{j}{mk} \leq x <\frac{i}{m}+\frac{j+1}{mk}.
\]
Put $\bar{f_1}(1) = f_0(1)$.
Now choose a finite system $\iI$ of pairwise disjoint open 
intervals covering $1$ and all
numbers of the form $\frac{i}{m}+\frac{j}{mk}$ such that
\beq\lab{ii}
\sum_{I\in\iI} |I|^{1-\al+\frac{1}{N}} < \frac{1}{2M}.
\eeq
We can clearly choose a continuous $f_1\in\C$ which is linear on each
$I\in\iI$, agrees with $\bar{f_1}$ outside $\cup\iI$ and satisfies 
$f_1(0)=f_0(0)$ and $f_1(1)=f_0(1)$. 
Denote the supremum-norm of a not necessarily continuous function $f$ by 
$||f||$. One can easily check that $||f_0-\bar{f_1}|| \leq \frac{2}{5} r_0$
and that $||\bar{f_1}-f_1|| \leq \frac{2}{5} r_0$, hence 
$||f_0-f_1|| \leq \frac{4}{5} r_0$.
So if we put $r_1=\frac{r_0}{20k} \leq \frac{r_0}{5}$ then 
$B(f_1,r_1) \subset B(f_0,r_0)$.
Now we claim that the inequality
\beq\lab{gyok1}
\left(\frac{2}{mk}\right)^\al < \frac{r_0}{10k}
\eeq
implies that for 
every $f\in B(f_1,r_1)$, every H\"older continuous function $g$
of exponent $\al$ and constant $1$ and every 
fixed $0\leq i <m$, there is at most one
$0\leq j <k$ such that $\left(\{f=g\}\cap 
\left(\frac{i}{m}+\frac{j}{mk},\frac{i}{m}+\frac{j+1}{mk}\right)\right) 
\setminus \cup\iI \neq \emptyset$.
Indeed, by the concavity of the function $x^\al$ it is enough to check 
that for this fixed $i$ the graph of $f$ and $g$ cannot meet over two 
consecutive intervals of the $k$ intervals of length $\frac{1}{mk}$, 
but this is clear from the fact 
that the value of $\bar{f_1}$ `jumps' by $\frac{r_0}{5k}$, moreover 
$\frac{r_0}{5k}-2r_1=\frac{r_0}{10k}$, and from (\ref{gyok1}).

This means that $\{f=g\}$ can be covered by the elements of 
$\iI$ and by $m$ intervals of 
length $\frac{1}{mk}$. Before we use this fact to estimate 
$\iH^{1-\al+\frac{1}{N}}_\infty (\{f=g\})$, we need some more preparations.

By rearranging (\ref{gyok1}) we obtain
\[
k < \left(\frac{r_0}{10 \cdot 2^\al}\right)^\frac{1}{1-\al} m^\frac{\al}{1-\al},
\]
and we will also make sure that even the following holds.
\beq\lab{gyok2}
\frac{1}{2}\left(\frac{r_0}{10 \cdot 2^\al}\right)^\frac{1}{1-\al} 
m^\frac{\al}{1-\al} <
k < \left(\frac{r_0}{10 \cdot 2^\al}\right)^\frac{1}{1-\al} m^\frac{\al}{1-\al}
\eeq

The last condition we need for the estimation of the Hausdorff capacity is 
that 
\beq\lab{pre}
m\left({\frac{1}{mk}}\right)^{1-\al+\frac{1}{N}} < \frac{1}{2M}.
\eeq
Indeed, together with (\ref{ii}) this implies that
$\iH^{1-\al+\frac{1}{N}}_\infty (\{f=g\}) < \frac{1}{M}$ for every 
$f\in B(f_1,r_1)$ and every H\"older continuous function $g$
of exponent $\al$ and constant $1$. 
Using the left
hand side of (\ref{gyok2}) we obtain
\beq\lab{pre2}
m\left({\frac{1}{mk}}\right)^{1-\al+\frac{1}{N}} <
\left(\frac{1}{2}\left(\frac{r_0}{10 \cdot 2^\al}\right)^\frac{1}{1-\al}\right)
^{-(1-\al+\frac{1}{N})} 
m^{1-(1-\al+\frac{1}{N})-\frac{\al(1-\al+\frac{1}{N})}{1-\al}},
\eeq
and as the exponent
\[
1-(1-\al+\frac{1}{N})-\frac{\al(1-\al+\frac{1}{N})}{1-\al} = 
-\frac{1}{(1-\al) N}
< 0,
\]
we can choose a large enough $m$ such that (\ref{pre}) holds, and then
we can fix $k$ according to (\ref{gyok2}), so the proof is complete.
\ep

However, Theorem \ref{Holder} is an upper estimate of the dimension, 
while in view of the question of isomorphisms 
we are more interested in lower estimates. Unfortunately we cannot
prove any.

\bq
Fix $0<\al<1$. Is the typical continuous function H\"older
continuous of exponent $\al$ on a set of positive Hausdorff dimension?
Or on a set of Hausdorff dimension $1-\al$?
\eq

Motivated by the previous problem we can formulate some other natural and
interesting questions. See
\cite{BH} and \cite{HL} for similar results. E.g.~in \cite{HL} it is shown
that a typical continuous function agrees with a monotone function only on a
set of dimension 0. It is therefore natural to guess that the same holds with 
functions of bounded variation, as they are the
differences of monotone functions. However, we were able to prove only the
following.

\bt\lab{BV}
A typical continuous function does not agree with any function of 
bounded variation on any set of Hausdorff dimension larger than $\frac{1}{2}$.
\et

\bp
The proof is similar to that of Theorem \ref{Holder}. 
As usual, total variation of a function $g$ is denoted by 
\[
Var(g) = \sup\left\{\sum_{i=1}^n |g(x_{i})-g(x_{i-1})| : n\in\NN,\  
0=x_0<x_1<\dots<x_n=1 \right\}.
\]
It is sufficient to prove that for all integers $N,M>0$ the set
\beq\lab{var}
\left\{f\in\C : \iH^{\frac{1}{2}+\frac{1}{N}}_\infty (\{f=g\})<\frac{1}{M} 
\textrm{ for every } g \textrm{ with } Var(g) \leq 1\right\}
\eeq
contains a dense open set. Given $f_0$ and $r_0$ define 
$\bar{f_1}$ as in Theorem \ref{Holder} with $m$ and $k$ yet unspecified. 
(In fact, we will choose $m=k$ so that 
$|x-y|<\frac{2}{m}$ implies $|f_0(x)-f_0(y)|<\frac{r_0}{5}$ and so that
the right hand side of (\ref{bv}) below is smaller than $\frac{1}{2M}$.)
 
Choose $\iI$ such that
\beq\lab{ii2}
\sum_{I\in\iI} |I|^{\frac{1}{2}+\frac{1}{N}} < \frac{1}{2M},
\eeq
and define $f_1$ and $r_1$ as before. 

Now for every $f\in B(f_1,r_1)$, $g$ as in (\ref{var}), and fixed 
$0 \leq i <m$, denote $l_i$ the number of intervals over which the graph of $f$ 
and $g$ meet outside $\cup\iI$; that is,
\[l_i = \#\left\{ j : \left(\{f=g\} \cap 
\left(\frac{i}{m}+\frac{j}{mk},\frac{i}{m}+\frac{j+1}{mk}\right)\right)
\setminus\cup\iI\neq \emptyset \right\}.
\]
As the `jumps' of $\bar{f_1}$ are 
of height $\frac{r_0}{5k}$, moreover $\frac{r_0}{5k}-2r_1=\frac{r_0}{10k}$,
it is easily seen that $Var(g) \geq 
\sum_{i=1}^{m-1} \left( \frac{r_0}{10k}(l_i-1) \right) = 
\frac{r_0}{10k} \left(\sum_{i=1}^{m-1} l_i - m\right)$. Thus 
$Var(g)\leq 1$ 
implies that $\sum_{i=1}^{m-1} l_i$; that is, the number intervals of length $\frac{1}{mk}$ needed to cover $\{f=g\} \setminus \cup\iI$ is at most 
$m+\frac{10k}{r_0}$.

Choose $k=m$. Then $\iH^{\frac{1}{2}+\frac{1}{N}}_\infty 
(\{f=g\})$ can be estimated by (\ref{ii2}) and
\beq\lab{bv}
\left(m+\frac{10m}{r_0}\right)\left(\frac{1}{m^2}\right)^{\frac{1}{2}+\frac{1}{N}} = 
\left(\frac{10}{r_0}+1\right) m^{-\frac{2}{N}},
\eeq
which is smaller than $\frac{1}{2M}$ if $m$ is large enough.
\ep

\bq
Does a typical continuous function agree with some function of 
bounded variation on some set of Hausdorff dimension $\frac{1}{2}$? 
Or on a set of positive dimension?
\eq

\bigskip

\noindent
\textsc{R\'enyi Alfr\'ed Institute, Re\'altanoda u.~13-15.~Budapest
1053, Hungary}

\textit{Email address}: \verb+emarci@renyi.hu+


\end{document}